# The Basel Problem

## Abstract

Because of its relation to the distribution of prime numbers, the Riemann zeta function ζ (s) is one of the most important functions in mathematics. The zeta function is defined by the following formula for any complex number s with the real component greater than 1. $\zeta(s) = \sum_{n=1}^{\infty} 1/n^s$ .Taking s=2, we see that ζ(2) is equal to the sum of the squares of reciprocals of all positive integers. This leads to the famous problem by Basel in mathematical analysis with important relevance to number theory, solved by Leonhard Euler in 1734. In this paper we discuss some of the notable proofs given by mathematicians to the basal problem. Most of the theorems are very well known whereas some can be found as proofs of problems present in textbooks.

## 1. Introduction

The Basel problem was first introduced in 1644 by Pietro Mengoli, Italian mathematician and clergyman (1626–1686) who is known for his Known (nowadays) for work in infinite series. The Problem remained open for 90 years, until by solving it, Euler gave his first proof in 1734. The Basel problem asks for the sum of the series:

$$\sum_{n=1}^{\infty} \frac{1}{n^2} = \frac{1}{1^2} + \frac{1}{2^2} + \frac{1}{3^2} + \cdots \tag{1}$$

It was Euler who found the exact sum to be $\frac{\pi^2}{6}$. In this paper we discuss some of the notable proofs given by mathematicians to the basal problem. Most of the theorems are very well known whereas some can be found as proofs of problems present in textbooks. The R.H.S. of (1) can be written as

$$\sum_{n=0}^{\infty} \frac{1}{(2n+1)^2} = \sum_{n=1}^{\infty} \frac{1}{n^2} - \sum_{n=1}^{\infty} \frac{1}{(2n)^2} = \left(1 - \frac{1}{4}\right) \sum_{n=1}^{\infty} \frac{1}{n^2} = \frac{3}{4} \sum_{n=1}^{\infty} \frac{1}{n^2} \tag{2}$$

**Proof 1:**

Let us first take a proof that is due to Calabi, Beukers and Kock[1] and also mentioned in Robin Chapman[14].
They first noticed the following equality:

$$\sum_{k=0}^{\infty} \frac{1}{(2k+1)^2} = \int_0^1 \int_0^1 \frac{dxdy}{1-x^2y^2} \tag{3}$$

After that they made the substitution $u = tan^{-1} x\sqrt{\frac{1-y^2}{1-x^2}}$ $and$ $v = tan^{-1} y\sqrt{\frac{1-x^2}{1-y^2}}$. We now get the transformed variables after substitution as $x = \frac{\sin u}{\cos v}$ $and$ $y = \frac{\sin v}{\cos u}$ . Now notice the following important equation which connects $\zeta(2)$ and the left hand side of the equation (3).

$$\frac{3}{4}\zeta(2) = \sum_{n=1}^{\infty} \frac{1}{n^2} - \sum_{m=1}^{\infty} \frac{1}{2m^2} = \sum_{k=0}^{\infty} \frac{1}{(2k+1)^2} \tag{4}$$

The Right hand side of the equation (3) can be rewritten after the substitution as

$$\int_0^1 \int_0^1 \frac{dxdy}{1-x^2y^2} = \int \int_A \frac{dudv}{1-u^2v^2} J \qquad here\ A = \{(u,v): u > 0, v > 0, u + v < \frac{\pi}{2}\} \tag{5}$$

Since the $J$ the Jacobian matrix will be

$$J = \frac{\partial(x,y)}{\partial(u,v)} = \det \begin{vmatrix} \cos u / \cos v & \sin u \sin v / \cos^2 v \\ \sin u \sin v / \cos^2 u & \cos v / \cos u \end{vmatrix} = 1 - u^2v^2 \tag{6}$$

Now as the area of $A = \left(\frac{1}{2}\right) * \left(\frac{\pi}{2}\right) * \left(\frac{\pi}{2}\right) = \frac{\pi^2}{8}$. We get from equation (1),(2) and (3) that $\zeta(2) = \frac{\pi^2}{6}$.

**Proof 2:**

We now derive a proof that comes from a note by Boo Rim Choe[2] in the American Mathematical Monthly in 1987 and also mentioned in Robin Chapman[14]. We start by taking a note of the power series expansion of the inverse sine function.

$$\sin^{-1}(x) = \sum_{n=0}^{\infty} \frac{1.3\ldots(2n-1)}{2.4\ldots(2n)} \frac{x^{2n+1}}{2n+1} \qquad for \qquad |x| \leq 1 \tag{7}$$

Now let us put $t = \sin^{-1}(x)$. After the substitution we get

$$t = \sum_{n=0}^{\infty} \frac{1.3\ldots(2n-1)}{2.4\ldots(2n)} \frac{\sin(t)^{2n+1}}{2n+1} \qquad for \qquad |t| \leq \pi/2 \tag{8}$$

To proceed further we first need to find the integration from $0\ to\ \pi/2$ of $\sin(t)^{2n+1}$ or $I_{2n+1} = \int_0^{\pi/2} (\sin(x))^{2n+1}\, dx$, we need to use the recursion formula which can be expressed as:

$$I_{2n+1} = I_{2n-1} - \int_0^{\pi/2} [(\sin(x))^{2n-1} \cos(x)]\, cosx.dx \tag{9}$$

Now using integration by parts now we can derive that

$$\frac{2n+1}{2n} I_{2n+1} = I_{2n-1} \tag{10a}$$

Therefore by recursion we can write that

$$I_{2n+1} = \frac{2n}{2n+1} I_{2n-1} = \frac{2n.2(n-1)}{(2n+1)(2n-1)} I_{2n-3} = \frac{2.4\ldots(2n)}{3.5\ldots(2n+1)} \tag{10b}$$

We therefore get that

$$I_{2n+1} = \int_0^{\pi/2} (\sin(x))^{2n+1}\, dx = \frac{2.4\ldots(2n)}{3.5\ldots(2n+1)} \tag{11}$$

Now if we use equation (11) and (8) we get

$$\frac{\pi^2}{8} = \int_0^{\pi/2} t.dt = \sum_{n=0}^{\infty} \frac{1}{(2n+1)^2} \tag{12}$$

The right hand side of the equality (12) is $(3/4)\zeta(2)$. Therefore we get that $\zeta(2) = \frac{\pi^2}{6}$.

**Proof 3:**

We now take a proof from an article in the Mathematical Intelligencer by Apostol[3] in 1983 and also mentioned in Robin Chapman[14]. We first start noting the following equation:

$$\frac{1}{j^2} = \int_0^1 \int_0^1 x^{j-1} y^{j-1} dxdy \tag{13}$$

It is also to note that we get by the monotone convergence theorem

$$\sum_{k=1}^{\infty} \frac{1}{k^2} = \int_0^1 \int_0^1 \left( \sum_{k=1}^{\infty} (xy)^{(n-1)} \right) . dxdy = \int_0^1 \int_0^1 \frac{dxdy}{1-xy} \tag{14}$$

Let us now introduce change of variables by changing the variables $u = (x + y)/2)$ $and$ $v = (x - y)/2)$). After the substitution of the variables we get the relationship between changed variables and actual variables as $x = u - v$ and $= u + v$. Now notice that the Left hand side of the equation (14) is $\zeta(2)$. We therefore can rewrite equation (14) in terms of the changed variables as:

$$\zeta(2) = 2 \iint_S \frac{dudv}{1-u^2+v^2} \tag{15}$$

Here S is the square with vertices (0, 0), (1/2, −1/2), (1, 0) and (1/2, 1/2). Now because of the symmetry present (symmetrical about the line y=0). We can find the value by only doubling the integration of the upper half of the square which can be itself be divided into two triangles with coordinates $((0,0), \left(\frac{1}{2}, \frac{1}{2}\right), \left(\frac{1}{2}, 0\right))$ and $(\left(\frac{1}{2}, \frac{1}{2}\right), \left(\frac{1}{2}, 0\right), (1,0))$ respectively. We can therefore express $\zeta(2)$ as

$$\zeta(2) = 4 \int_0^{\frac{1}{2}} \int_0^u \frac{dudv}{1-u^2+v^2} + 4 \int_{\frac{1}{2}}^1 \int_0^{1-u} \frac{dudv}{1-u^2+v^2}$$

$$= 4 \int_0^{\frac{1}{2}} \frac{1}{\sqrt{1-u^2}} tan^{-1} \left( \frac{u}{\sqrt{1-u^2}} \right) . du + 4 \int_{1/2}^1 \frac{1}{\sqrt{1-u^2}} tan^{-1} \left( \frac{1-u}{\sqrt{1-u^2}} \right) . du \tag{16}$$

Now we notice that $tan^{-1} \left( \frac{u}{\sqrt{1-u^2}} \right) = sin^{-1}(u)$. Let us this identity by again changing the variables by substituting $\theta = tan^{-1}(\frac{1-u}{\sqrt{1-u^2}})$. After substitution we get

$$tan^2\theta = (1-u)/(1+u) \text{ and } sec^2\theta = \frac{2}{1+u}. \tag{17}$$

It follows that $u = 2cos^2\theta - 1 = cos2\theta$ and therefore we get that $\theta = \left(\frac{1}{2}\right) cos^{-1}u = \frac{\pi}{4} - \frac{1}{2} sin^{-1}u$. We can hence rewrite equation (16) in terms of the changed variables as

$$\zeta(2) = 4 \int_0^{\frac{1}{2}} \frac{sin^{-1}u}{\sqrt{1-u^2}} du + 4 \int_{\frac{1}{2}}^1 \frac{1}{\sqrt{1-u^2}} \left( \frac{\pi}{4} - \frac{sin^{-1}u}{2} \right) . du = [2(sin^{-1}u)^2]_0^{1/2} + [\pi sin^{-1}u - (sin^{-1}u)^2]_{1/2}^1 = \frac{\pi^2}{18} + \frac{\pi^2}{2} - \frac{\pi^2}{4} - \frac{\pi^2}{6} + \frac{\pi^2}{36} = \frac{\pi^2}{6}. \tag{18}$$

**Proof 4:**

We now present some textbook proofs, found in many books on Fourier analysis to find the value of $\zeta(2)$. In this proof we are going to discuss next we will use the $L^2$ completeness of the trigonometric functions. Suppose that $\boldsymbol{H}$ is a Hilbert space with inner product. Let $(e_n)$ be an ortho-normal basis of H such that

$$\langle e_m, e_n \rangle = \begin{cases} 1 & if \quad m = n \\ 0 & if \; m \neq n \end{cases} \tag{19}$$

Then we know that by Parseval's identity, we can write for every $g \in \boldsymbol{H}$,

$$\langle g,g\rangle = ||g||^2 = \sum_n |\langle g,e_n\rangle|^2 \tag{20}$$

Let the Hilbert Space $\boldsymbol{H}$ be $L^2[0,1]$. Let us now apply this theorem $g(x) = x$. The Left hand side of the equation (20) will be $\langle g,g\rangle = \frac{1}{3}$ and $\langle g,e_0\rangle = \frac{1}{2}$ and $\langle g,e_n\rangle = \frac{1}{2\pi in}$ for $n \neq 0$. Therefore we can re express equation (20) as:

$$\frac{1}{3} = \frac{1}{4} + \sum_{\boldsymbol{n\in Z,n\neq 0}} \frac{1}{4\pi^2 n^2} \tag{21}$$

Which therefore gives us $\zeta(2) = \pi^2/6$.

**Proof 5:**

Let us now see another similar way of proving the identity by fourier transform. We first take a function f which is continuous in [ 0, 1 ] and $f(0) = f(1)$. Let us first explore the criteria for point-wise convergence of a periodic function f. The criteria for point-wise convergence of a periodic function f are as follows:

- Fourier series converges uniformly if $f$ satisfies a Holder condition.
- Fourier series converges everywhere if $f$ is of bounded variation.
- Fourier series converges uniformly if $f$ is continuous and its Fourier coefficients are absolutely summable.

Then the Fourier series of $f$ converges to $f$ pointwise. As our function holds true in both the circumstances we can say that the Fourier series of f converges to f point wise. Let us apply this to the function with these properties let's say $f(x) = x(1 - x)$ which then gives us the following equation:

$$x(1 - x) = \frac{1}{6} - \sum_{n=1}^{\infty} \frac{\cos(2\pi nx)}{\pi^2 n^2} \tag{22}$$

If we put x=0 in the LHS of the equation (22) we get $(2) = \pi^2/6$. We can also get the identity if we put x=1/2 which gives us:

$$\frac{\pi^2}{12} = -\sum_{n=1}^{\infty} \frac{(-1)^n}{n^2} \tag{23}$$

Now as we know that

$$\zeta(2) = \sum_{n=1}^{\infty} 1/n^2 \tag{24}$$

Now if we subtract the equation (23) from the equation (24). We get

$$\zeta(2) - \frac{\pi^2}{12} = 2 * \left(\sum_{n=1}^{\infty} \frac{1}{(2n)^2}\right) = \frac{1}{2} * \sum_{n=1}^{\infty} 1/n^2 = \zeta(2)/2 \tag{25}$$

We therefore get $\zeta(2)/2 = \frac{\pi^2}{12}$ which gives us $\zeta(2) = \pi^2/6$.

This is mentioned in Robin Chapman[14].

**Proof 6:**

Let us now derive the original proof given by Euler. We use the infinite product

$$\sin(\pi x) = \pi x \prod_{n=1}^{\infty} \left(1 - \frac{x^2}{n^2}\right) \tag{26a}$$

for the sine function.

$$\sin \pi x = \pi x(1-x^2)\left(1-\frac{x^2}{4}\right)\left(1-\frac{x^2}{9}\right)\left(1-\frac{x^2}{16}\right) \dots \dots \quad (26b)$$

Now the RHS can be written as

$$RHS = \pi x + \pi x^3\left(1+\frac{1}{4}+\frac{1}{9}+\frac{1}{16}+\cdots\right)+\pi x^5\left(\frac{1}{1.4}+\frac{1}{1.9}+..+\frac{1}{4.9}+..\right)+.. \quad (27a)$$

Whereas the LHS can be expanded in power series as:

$$LHS = \pi x - \frac{(\pi x)^3}{3!} + \frac{(\pi x)^5}{5!} \quad (27b)$$

Now if we compare the coefficients of $x^3$ in the MacLaurin series of both the LHS and RHS sides we get $\zeta(2) = \pi^2/6$.

**Proof 7:**

The two earlier proofs using Fourier transform can be proved in a similar way but without using Fourier transform in the following manner. Consider the series:

$$f(t) = \sum_{n=1}^{\infty} \cos nt/n^2 \quad (28)$$

First of all notice that the function is uniformly convergent on the real line. Now we know that $\sin nt$ can also be written as $(e^{int}-e^{-int})/2i$ . Therefore we can write as follows:

$$\sum_{n=1}^{N} \sin nt = \sum_{n=1}^{N}(e^{int}-e^{-int})/2i = \frac{e^{it}-e^{i(N+1)t}}{2i(1-e^{it})} - \frac{e^{-it}-e^{-i(N+1)t}}{2i(1-e^{-it})}$$

$$= \frac{e^{it}-e^{i(N+1)t}}{2i(1-e^{it})} - \frac{1-e^{-iNt}}{2i(1-e^{it})} \quad (29)$$

And so this sum is bounded above in absolute value by $\frac{2}{|1-e^{it}|} = \frac{1}{\sin t/2}$. Hence these sums are uniformly bounded and by Dirichlet's test the sum $\sum_{n=1}^{N} \frac{\sin nt}{n}$ is uniformly convergent. It therefore follows that for $t \in (0, 2\pi)$ the following equation holds

$$f'(t) = -\sum_{n=1}^{\infty} \sin\frac{nt}{n} = -Im\left(\sum_{1}^{\infty}\frac{e^{int}}{n}\right) = Im(\log(1-e^{it})) = arg\,(1-e^{it}) = \frac{t-\pi}{2} \quad (30)$$

Now notice that

$$f(\pi)-f(0) = \int_0^{\pi} f'(t)dt = \int_0^{\pi}\frac{t-\pi}{2}dt = -\frac{\pi^2}{4} \quad (31)$$

But $f(0) = \zeta(2)$ and $f(\pi) = \sum_{n=1}^{\infty}\frac{(-1)^n}{n^2} = -\zeta(2)/2$. Hence $\zeta(2) = \pi^2/6$.

This is mentioned in Robin Chapman[14].

**Proof 8:**

Let us discuss about another textbook proof, found in many books on complex analysis. We use the calculus of residues. Let $g(x) = \pi z^{-2} cot\, \pi z$. It can be easily seen that $g$ has poles at precisely the integers where the value of $cot\, \pi z$ becomes infinite.

Now notice that at the pole at zero $g(x)$ has residue $-\pi^2/3$, and that at a non-zero integer has residue $1/n^2$. Let N be a natural number and let $C_N$ be the square contour with vertices $(\pm 1 \pm i)(N + 1/2)$. Now we know that by the calculus of residues the following identity holds:

$$\frac{-\pi^2}{3} + 2\sum_{n=1}^{N} \frac{1}{n^2} = \frac{1}{2\pi i}\int_{C_N} g(z).dz = I_n \tag{32}$$

Now if $\pi z = x + iy$ the following equation holds:

$$|\cot \pi z|^2 = \frac{\cos^2 x + \sinh^2 x}{\sin^2 x + \sinh^2 x} \tag{33}$$

Now if z lies on the vertical edges of $C_N$ then the following equation holds

$$|\cot \pi z|^2 = \frac{\sinh^2 x}{1 + \sinh^2 x} < 1 \tag{34}$$

and if z lies on the horizontal edges of $C_N$ then the following equation holds

$$|\cot \pi z|^2 \le \frac{1 + \sinh^2 \pi(N + \frac{1}{2})}{\sinh^2 \pi(N + \frac{1}{2})} = \coth^2 \pi(N + \frac{1}{2}) \le \coth^2(\pi/2) \tag{35}$$

Hence $|\cot \pi z| \le K = \coth(\pi/2)$ on $C_N$, and so $|f(z)| \le \frac{\pi K}{(N+\frac{1}{2})^2}$ on $C_N$. This estimate shows that the following equation holds:

$$|I_n| \le \frac{8\pi K(N+\frac{1}{2})}{2\pi(N+\frac{1}{2})^2} \tag{36}$$

Notice that as $I_n \to 0$ as $N \to \infty$. Therefore we get $\zeta(2) = \pi^2/6$.

This is mentioned in Robin Chapman[14].

**<u>Proof 9:</u>**

Let us discuss another proof that can be found as an exercise in Apostol[4] and also mentioned in Robin Chapman[14].We first start by taking the note of two important results

If $0 < x < \pi/2$ then the following 2 inequality holds:

- $\sin x < x < \tan x$ and so
- $\cot^2 x < x^{-2} < 1 + \cot^2 x$

We can therefore write that If n and N are natural numbers with $1 \le n \le N$, then

$$\cot^2 \frac{n\pi}{2N+1} < \frac{(2n+1)^2}{n^2\pi^2} < 1 + \cot^2 \frac{n\pi}{2N+1} \tag{37}$$

Using the first inequality of equation (37) we get

$$\frac{\pi^2}{(2n+1)^2}\sum_{n=1}^{N} \cot^2 \frac{n\pi}{2N+1} < \sum_{n=1}^{N} \frac{1}{n^2} < \frac{N\pi^2}{(2n+1)^2} + \frac{\pi^2}{(2n+1)^2}\sum_{n=1}^{N} \cot^2 \frac{n\pi}{2N+1} \tag{38}$$

Now if $A_n = \sum_{n=1}^{N} \cot^2 \frac{n\pi}{2N+1}$, we can prove the identity $\zeta(2) = \frac{\pi^2}{6}$ by simply showing that $\lim_{N\to\infty} \frac{A_n}{N^2} = \frac{2}{3}$. Now If $1 \le n \le N$ and $\theta = \frac{n\pi}{2N+1}$, then $sin(2N+1)\theta = 0$.

Now by De-moivre's theorem we know that $sin(2N+1)\theta$ is the imaginary part of $(\cos\theta + i\sin\theta)^{2N+1}$, and so

$$\frac{sin(2N+1)\theta}{sin^{2N+1}\theta} = \frac{1}{sin^{2N+1}\theta}\sum_{k=0}^{N}(-1)^k C_{2N-k}^{2N+1} cos^{2(N-k)}\theta\, sin^{2k+1}\theta = \sum_{k=0}^{N}(-1)^k C_{2N-k}^{2N+1} cot^{2(N-k)}\theta \tag{39}$$

This can be written in the form

$$f(x) = (2N+1)x^N - C_3^{2N+1}x^{N-1} + \cdots. \tag{40}$$

Where $f(cot^2\theta) = \sum_{k=0}^{N}(-1)^k C_{2N-k}^{2N+1} cot^{2(N-k)}\theta$.

Hence the roots of $f(x) = 0$ are $cot^2\,(n\pi/(2N+1))$ where $1 \le n \le N$ and so $A_n = N(2N-1)/3$. Thus $\frac{A_n}{N^2} \to \frac{2}{3}$ , as required.

**Proof 10:**

This comes from a note by Kortram[5] and also mentioned in Robin Chapman[14]. Given an odd integer n = 2m + 1 a very important identity is that $sin\, nx = F_n(sin\, x)$ where $F_n$ is a polynomial of degree $n$. An important property to notice is that the zeros of $F_n(y)$ are the values $sin(j\pi/n)$ $(-m \le j \le m)$ and $\lim_{y\to 0}\left(\frac{F_n(y)}{y}\right) = n$. We can then say that

$$F_n(y) = ny\prod_{j=1}^{m}\left(1 - \frac{y^2}{sin^2(\frac{j\pi}{n})}\right) \tag{41}$$

and therefore by our definition the following equality holds

$$\sin nx = nsinx\prod_{j=1}^{m}\left(1 - \frac{sin^2 x}{sin^2(\frac{j\pi}{n})}\right) \tag{42}$$

Now let us compare the coefficients of $x^3$ in the MacLaurin expansion of both sides LHS and RHS. Now the MacLaurin expansion of $\sin nx$ gives us:

$$\sin nx = nx - \frac{n^3x^3}{3!} + \frac{n^5x^5}{5!} + .. \tag{43}$$

And the coefficients of $x^3$ in RHS will be $-\frac{n}{6} - n\sum_{j=1}^{m}\frac{1}{sin^2(\frac{j\pi}{n})}$. Therefore we can write by equation (42) and (43) that$-\frac{n^3}{6} = -\frac{n}{6} - n\sum_{j=1}^{m}\frac{1}{sin^2(\frac{j\pi}{n})}$ .We can therefore say that using the earlier equations (43) the following identity holds:

$$\frac{1}{6n^2} = \frac{1}{6} - \sum_{j=1}^{m}\frac{1}{n^2 sin^2(\frac{j\pi}{n})} \tag{44}$$

Let us fix an integer $M$ and let $m > M$ without any loss of generality. Then we can say that

$$\frac{1}{6n^2} + \sum_{j=M+1}^{m}\frac{1}{n^2 sin^2(\frac{j\pi}{n})} = \frac{1}{6} - \sum_{j=1}^{M}\frac{1}{n^2 sin^2(\frac{j\pi}{n})} \tag{45}$$

and using the inequality $sin\ x\ > \frac{2}{\pi}x$ for $0\ <\ x\ < \frac{\pi}{2}$, we get

$$0 < \frac{1}{6} - \sum_{j=1}^{M} \frac{1}{n^2 sin^2(\frac{j\pi}{n})} < \frac{1}{6n^2} + \sum_{j=M+1}^{m} \frac{1}{4j^2} \tag{46}$$

Now notice that if we let n towards infinity we get

$$0 < \frac{1}{6} - \sum_{j=1}^{M} \frac{1}{\pi^2 j^2} < \sum_{j=M+1}^{\infty} \frac{1}{4j^2} \tag{47}$$

Hence we get easily get $\sum_{j=1}^{\infty} \frac{1}{\pi^2 j^2} = \frac{1}{6} \rightarrow \sum_{j=1}^{\infty} \frac{1}{j^2} = \frac{\pi^2}{6}$ as wanted.

**Proof 11:**

This proof again comes from Euler which uses the Weierstrass Factorization Theorem which is a direct generalization of the Fundamental Theorem of Algebra. We first state that theorem.

Let $f$ be an entire function and let $\{ a_n\}$ be the nonzero zeros of $f$. Suppose f has a zero at $z\ =\ 0$ of order $m\ \geq\ 0$ (where order 0 means $f(0)\ \neq\ 0$). Then there exist a function $g$ and a sequence of integers $\{p_n\}$ such that the following equation holds:

$$f(z) = z^m \exp\ (g(z)) \prod_{n=1}^{\infty} E_{p_n}(z/a_n)$$

Where $E_n(y) = (1-y) \qquad if \qquad n = 0;$

$E_n(y) = (1-y)\exp\ (y + \frac{y^2}{2} + \cdots + \frac{y^n}{n}) \qquad if \qquad n = 1,2,3..;$

This. Let us assume $f\ =\ sin(\pi x)$, the sequence $p_n = 1$ and the function $g(z)\ =\ log(\pi z)$. Therefore we can write as follows:

$$\sin t = t\left(1-\frac{t}{\pi}\right)\left(1+\frac{t}{\pi}\right)\left(1-\frac{t}{2\pi}\right)(1+\frac{t}{2\pi}) \tag{48}$$

Since $sin\ t$ have roots precisely at $t\ \in\ Z$. Now if we substitute $t = \pi y$ we get the following equation:

$$\sin(\pi y) = \pi y(1-y)(1+y)\left(1-\frac{y}{2}\right)\left(1+\frac{y}{2}\right)\ldots = \pi y(1-y^2)\left(1-\frac{y^2}{4}\right)\left(1-\frac{y^2}{9}\right).. \tag{49}$$

Taking logarithm of both sides of the equality

$$\ln(\sin(\pi y)) = ln\pi + lny + \ln(1-y^2) + [ln(4-y^2) - ln4] + \cdots \tag{50}$$

Now if we differentiate both sides of the equality with the variable $y$ we get the following equation:

$$\pi\cos(\pi y) * \left(\frac{1}{\sin(\pi y)}\right) = \frac{1}{y} - \frac{2y}{1-y^2} - \frac{2y}{4-y^2} - \frac{2y}{9-y^2} - \cdots \tag{51}$$

This gives us

$$\frac{1}{y} + \frac{1}{1-y^2} + \frac{1}{4-y^2} + \frac{1}{9-y^2} + \cdots. = \frac{1}{2y^2} - \pi\cos(\pi y) * \left(\frac{1}{2y\sin(\pi y)}\right) \tag{52}$$

Now we put $y = -ix$, we can write the above equation (52) as follows:

$$\frac{1}{1+x^2}+\frac{1}{4+x^2}+\frac{1}{9+x^2}+\cdots=-\frac{1}{2x^2}+\pi\cos(-i\pi x)*\left(\frac{1}{2\text{ixsin}(-i\pi x)}\right) \tag{53}$$

Use Euler's Formula we can write that

$$\frac{cos(z)}{sin(z)}=\frac{\frac{1}{2}(e^{iz}+e^{-iz})}{\frac{1}{2}(e^{iz}-e^{-iz})}=\frac{i(e^{2iz}+1)}{e^{2iz}-1} \tag{54}$$

$$\frac{\pi cos(-i\pi x)}{2ixsin(-i\pi x)}=\frac{\pi}{2ix}.\frac{i(e^{2\pi x}+1)}{e^{2\pi x}-1}=\frac{\pi}{2x}.\frac{e^{2\pi x}+1}{e^{2\pi x}-1}=\frac{\pi}{2x}+\frac{\pi}{x(e^{2\pi x}-1)} \tag{55}$$

Now if we use equation (54) and (55) then we can write as follows:

$$\frac{1}{1+x^2}+\frac{1}{4+x^2}+\frac{1}{9+x^2}+\cdots=-\frac{1}{2x^2}+\pi\cos(-i\pi x)*\left(\frac{1}{2\text{ixsin}(-i\pi x)}\right)=-\frac{1}{2x^2}+\frac{\pi}{2x}+\frac{\pi}{x(e^{2\pi x}-1)}=\frac{e^{2\pi x}(\pi x-1)+\pi x+1}{2x^2(e^{2\pi x}-1)} \tag{56}$$

LHS of the equation can be easily be transformed into $\sum_{j=1}^{\infty}\frac{1}{j^2}$ by substituting$x=0$. Notice that the form is $\frac{0}{0}$ in the Right hand side of the equation (56) .Therefore we can use L'hospital's rule which yields $\zeta(2)=\pi^2/6$.

**Proof 12:**

This is proof is due to Luigi Pace, Dept of Econ & Stats at Udine, Italy[6] .It was inspired by a 2003 note that solves the problem using a double integral on $R_+^2$ via Fubini's Theorem; a result that gives conditions under which it is possible to compute a double integral by using an iterated integral.

We know define two variables $X_1, X_2\ :\ R\rightarrow\ R_+$ random variables having a probability density function: $\rho_{X_i}: R_+\rightarrow\ [0,1]$ ]. Let us now define another variable $Y\ =\ X_1/X_2$ . We claim that Probability density function for $Y$ is $\rho_Y=\int_0^\infty t\,\rho_{X_1}(t\mu)\rho_{X_2}(t).\,dt$. To prove that joint density:

$$\Pr(a\le Y\le b)=\int_0^\infty\int_{at_2}^{bt_2}\rho_{X_1}(t_1)\rho_{X_2}(t_2).\,dt_1\,dt_2=\int_0^\infty\int_a^b t_2\rho_{X_1}(t_2u)\rho_{X_2}(t_2).\,du\,dt_2=\int_a^b\int_0^\infty t_2\rho_{X_1}(t_2u)\rho_{X_2}(t_2).\,du\,dt_2 \tag{57}$$

For the proof let us assign half-Cauchy distribution to $X_1, X_2$independently: $\rho_{X_i}(t)=\frac{2}{\pi(1+t)^2}$. Now we shall use this in the formula for $\rho_Y(u)$ obtained above:

$$\rho_Y(u)=\frac{4}{\pi^2}\int_0^\infty\frac{t}{1+t^2u^2}.\frac{1}{1+t^2}dt=\frac{2}{\pi^2(u^2-1)}\left[\ln\left(\frac{1+t^2u^2}{1+t^2}\right)\right]_0^\infty=\frac{4}{\pi^2}\frac{\ln u}{u^2-1} \tag{58}$$

Integrating $\rho_Y(u)$ from 0 to 1 we get $Pr(0\ \le\ Y\ \le\ 1)$. Therefore we can say that

$$Pr(0\ \le\ Y\ \le\ 1)=\int_0^1\rho_Y(u).\,du=\int_0^1\frac{4}{\pi^2}\frac{\ln u}{u^2-1}.\,du \tag{59}$$

But note that $Pr(0\ \le\ Y\ \le\ 1)\ =\ ½$ , Therefore we can say using equations (59) that

$$\int_0^1\frac{\ln u}{u^2-1}.\,du=\frac{\pi^2}{8}. \tag{60}$$

Let us simplify the integral by using $\frac{1}{1-u^2}=1+u^2+u^4+\cdots$. We can therefore get as follows:

$$\frac{\pi^2}{8} = -\int_0^1 \frac{\ln u}{1-u^2}.du = \sum_{n=0}^{\infty}\int_0^1 \frac{\ln u}{u^{2n}}.du = \sum_{n=0}^{\infty}\frac{1}{(2n+1)^2} = \frac{3}{4}\zeta(2) \tag{61}$$

Which gives us as required $\zeta(2) = \pi^2/6$.

**Proof 13:**

This proof is due to Matsuoka[7] and also mentioned in Robin Chapman[14].Consider the two integrals

$$I_n = \int_0^{\frac{\pi}{2}} cos^{2n}x.dx \ and \ J_n = \int_0^{\frac{\pi}{2}} x^2 cos^{2n}x.dx \tag{62}$$

Now as we discussed earlier

$$I_n = \frac{1.3.5....(2n-1)}{2.4.6...2n}\frac{\pi}{2} = \frac{(2n)!\pi}{4^n n!^2 2} \tag{63}$$

If $n > 0$ then we can say that if we do integration by parts gives

$$I_n = [xcos^{2n}x]_0^{\frac{\pi}{2}} + 2n\int_0^{\frac{\pi}{2}} xsinxcos^{2n-1}x = n[x^2 sinxcos^{2n-1}x]_0^{\frac{\pi}{2}} - n\int_0^{\frac{\pi}{2}} x^2(cos^{2n}x - (2n-1)sin^2xcos^{2n-2}x).dx = n(2n-1)J_{n-1} - 2n^2 J_n \tag{64}$$

we can therefore write using equation (64) and (63) that

$$\frac{(2n)!\pi}{4^n n!^2 2} = n(2n-1)J_{n-1} - 2n^2 J_n \tag{65}$$

Therefore we can rewrite expression (65) as:

$$\frac{\pi}{4n^2} = \frac{4^{n-1}[(n-1)!]^2}{(2n-2)!}J_{n-1} - \frac{(4^n n!^2)}{(2n)!}J_n \tag{66a}$$

Now summing up both the sides of the equation (66) we get:

$$\frac{\pi}{4}\sum_{n=1}^{N}\frac{1}{n^2} = \sum_{n=1}^{N}\frac{4^{n-1}[(n-1)!]^2}{(2(n-1))!}J_{n-1} - \sum_{n=1}^{N}\frac{(4^n n!^2)}{(2n)!}J_n = J_0 + \left[\frac{4^{2-1}[(2-1)!]^2}{(2(2-1))!}J_{2-1} - \frac{(4^1 1!^2)}{(2)!}J_1\right] + .. - \frac{(4^n N!^2)}{(2N)!}J_N \tag{66b}$$

All the middle terms present in the sequence gives a value a zero, and we are left with only the first term and the last term. Hence we can write that

$$\frac{\pi}{4}\sum_{n=1}^{N}\frac{1}{n^2} = J_0 - \frac{(4^n N!^2)}{(2N)!}J_N \tag{67}$$

Since $J_0 = \pi^3/24$ it sufficient to show that $\lim_{N\to\infty}\frac{(4^n N!^2)}{(2N)!} = 0$ to prove our identity.

Notice that the inequality $x < \frac{\pi}{2}sin\ x$ for $0 < x < \frac{\pi}{2}$ gives

$$J_N < \frac{\pi^2}{4}\int_0^{\frac{\pi}{2}} sin^2x.cos^{2n}x.x.dx = \frac{\pi^2}{4}\left(\int_0^{\frac{\pi}{2}} cos^{2n}x.dx - \int_0^{\frac{\pi}{2}} cos^{2n+2}x.dx\right) = \frac{\pi^2}{4}(I_N - I_{N+1}) = \frac{\frac{\pi^2}{8}.I_N}{N+1} \tag{68}$$

Thus we can write

$$0 < \frac{(4^n N!^2)}{(2N)!}J_N < \frac{\pi^3}{16(N+1)} \tag{69}$$

Hence proved that $\zeta(2) = \pi^2/6$.

**Proof 14:**

This proof is due to Stark[8] and also mentioned in Robin Chapman[14]. Consider a very important property of identity of the Fejer kernel. Before that we first give the definition of Fejer kernel. Fejér kernel is defined as

$$F_n(x) = \frac{1}{n}\sum_{k=0}^{n-1} D_k(x) \tag{70}$$

Where $D_k(x) = \sum_{s=-k}^{k} e^{isx}$ which is also known as the $k^{\text{th}}$ order Dirichlet Kernel. It can also be expressed as

$$F_n(x) = \sum_{|j|\leq n-1}(1-\frac{|j|}{n})e^{ijx} = \sum_{j=-n}^{n}(\mathrm{n}-|\mathrm{j}|)\mathrm{e}^{\mathrm{ijx}} \tag{71}$$

Using equation (70) and (71) we get

$$(\frac{\sin\frac{nx}{2}}{\sin\frac{x}{2}})^2 = \sum_{k=-n}^{n}(n-|k|)e^{ikx} = n + 2\sum_{k=1}^{n}(n-k)\cos kx \tag{72}$$

Hence we can write

$$\boldsymbol{I_n} = \int_{\mathbf{0}}^{\boldsymbol{\pi}} x(\frac{\sin\frac{nx}{2}}{\sin\frac{x}{2}})^2 . dx = \frac{n\pi^2}{2} + 2\sum_{k=1}^{n}(n-k)\int_0^\pi x\cos kx . dx = \frac{n\pi^2}{2} - 2\sum_{k=1}^{n}(n-k)\frac{1-(-1)^k}{k^2} = \frac{n\pi^2}{2} - 4n\sum_{1\leq k\leq n, k \text{ is a odd number}}\frac{1}{k^2} + 4\sum_{1\leq k\leq n, k \text{ is a odd number}}\frac{1}{k} \tag{73}$$

If we assume that $n = 2N$ with $N$ an integer then by equation (73) we can say that

$$\int_0^\pi \frac{x}{8N}(\sin Nx/\sin\frac{x}{2})^2 . dx = \frac{\pi^2}{8} - \sum_{r=0}^{N-1}\frac{1}{(2r+1)^2} + O(\frac{logN}{N}) \tag{74}$$

Now for $0 < x < \pi$ we can say that

$$sin\, x/2 > x/\pi \tag{75}$$

Therefore we can write the following equation:

$$\int_0^\pi \frac{x}{8N}(\sin Nx/\sin\frac{x}{2})^2 . dx < \frac{\pi^2}{8N}\int_0^\pi sin^2 Nx\frac{dx}{x} = \frac{\pi^2}{8N}\int_0^{N\pi} sin^2 y\frac{dy}{y} = O(logN/N) \tag{76}$$

Taking limits as $N \rightarrow \infty$ gives us the following equation:

$$\frac{\pi^2}{8} = \sum_{r=0}^{\infty}\frac{1}{(2r+1)^2} = \frac{3}{4}\zeta(2) \tag{77}$$

Hence proved that $\zeta(2) = \pi^2/6$.

**Proof 15:**

This is an exercise in Borwein & Borwein's Pi and the AGM[9] and also mentioned in Robin Chapman[14]. We carefully square Gregory's formula which states that

$$\frac{\pi}{4} = 1 - \frac{1}{3} + \frac{1}{5} - \frac{1}{7} + \cdots . = \sum_{n=0}^{\infty}\frac{(-1)^n}{2n+1} \tag{78}$$

Let us now denote $a_n = \sum_{n=-N}^{n=N} \frac{(-1)^n}{2n+1}$. Now notice that the sum of the positive powers will be same as sum of the negative powers. Since

$$\sum_{n=-\infty}^{-1} \frac{(-1)^n}{2n+1} = \sum_{n=1}^{\infty} \frac{(-1)^n}{-2n+1} = \sum_{l=0}^{\infty} \frac{(-1)^{l+1}}{-2(l+1)+1} = -\sum_{l=0}^{\infty} \frac{(-1)^l}{-2l-1} = \sum_{l=0}^{\infty} \frac{(-1)^l}{2l+1} \tag{79}$$

Therefore $\lim_{N\to\infty} a_n = \sum_{n=-\infty}^{n=\infty} \frac{(-1)^n}{2n+1} = 2 * \sum_{n=0}^{\infty} \frac{(-1)^n}{2n+1} = \pi/2$ . Let us denote $b_n = \sum_{n=-N}^{n=N} \frac{1}{(2n+1)^2}$. Now to prove the identity we require that $\lim_{N\to\infty} b_n = \frac{\pi^2}{4}$ , so we have to show that $\lim_{N\to\infty} (a_n{}^2 - b_n) = 0$. We now use the following identity:

$$\frac{1}{(2n+1)(2m+1)} = \frac{1}{2(m-n)}\left(\frac{1}{2n+1} - \frac{1}{2m+1}\right) \quad if\ m \neq n \tag{80}$$

Now if we square $a_n$ we get

$$a_n{}^2 = \sum_{n=-N}^{n=N} \frac{1}{(2n+1)^2} + \sum_{n=-N}^{N} \sum_{m=-N}^{N} \frac{1}{(2n+1)(2m+1)} \tag{81}$$

and therefore we can say that $(a_n{}^2 - b_n) = \sum_{n=-N}^{N} \sum_{m=-N}^{N} \frac{(-1)^{m+n}}{2(m-n)} * \left(\frac{1}{2n+1} - \frac{1}{2m+1}\right) = \sum_{n=-N}^{N} \frac{(-1)^N c_{n,N}}{2n+1}$ (82)

Here a point to note is that all the terms with zero denominators are omitted which comes when $m = n$ , $c_{n,N}$ referred above in the equation (82) is defined as

$$c_{n,N} = \sum_{n=-N}^{N} \frac{(-1)^m}{(m-n)} \tag{83}$$

It is easy to see that $c_{-n,N} = -c_{n,N}$ and so $c_{0,N} = 0$. If n > 0 then $c_{n,N} = \sum_{j=N-n+1}^{N+n} \frac{(-1)^j}{j}$ and so $|c_{n,N}| \leq 1/(N - n + 1)$. Thus we can write that

$$(a_n{}^2 - b_n) = \sum_{n=1}^{N}\left(\frac{1}{(2N-1)(N-n+1)} - \frac{1}{(2N+1)(N-n+1)}\right) = \sum_{n=1}^{N} \frac{1}{2N+1}\left(\frac{2}{2n-1} + \frac{1}{(N-n+1)}\right) + \sum_{n=1}^{N} \frac{1}{2N+3}\left(\frac{2}{2n+1} + \frac{1}{(N-n+1)}\right) \leq \frac{1}{2N+1}(2 + 4\log(2N+1) + 2 + 2\log(N+1)) \tag{84}$$

And so a $(a_n{}^2 - b_n) \to 0$ as $N \to \infty$ as required.

**Proof 16:**

This is an exercise in Hua's textbook on number theory[10] and also mentioned in Robin Chapman[14]. This depends on the formula for the number of representations of a positive integer as a sum of four squares. Let $r(n)$ be the number of quadruples $(x, y, z, t)$ of integers such that $= x^2 + y^2 + z^2 + t^2$ . Trivially $r(0) = 1$ and it is well known that

$$r(n) = 8 \sum_{m|n, m\ is\ not\ divisible\ by\ 4} m \tag{85}$$

for n > 0. Let $R(N) = \sum_{n=0}^{N} r(n)$ . Now notice that $R(N)$ is asymptotic to the volume of the 4D ball of radius $\sqrt{N}$, i.e., $R(N) \sim (\pi^2/2).N^2$ . But we know that

$$R(N) = 1 + 8\sum_{n=1}^{N} \sum_{m|n, m\ is\ not\ divisible\ by\ 4} m = 1 + 8\sum_{m \leq N, m\ is\ not\ divisible\ by\ 4} m \left\lfloor \frac{N}{m} \right\rfloor = 1 + 8(\theta(n) - 4\theta(N/4)) \tag{86}$$

Where $\theta(\mathrm{x}) = \sum_{\mathrm{m}\leq\mathrm{x}} \mathrm{m}\left\lfloor\frac{\mathrm{x}}{\mathrm{m}}\right\rfloor$. But we know that $\theta(\mathrm{x}) = \sum_{mr\leq x} m = \sum_{r\leq x}\sum_{m=1}^{\lfloor x/r\rfloor} m = \frac{1}{2}.\sum_{\mathrm{r}\leq\mathrm{x}}(\left\lfloor\frac{\mathrm{x}}{\mathrm{r}}\right\rfloor + \left\lfloor\frac{\mathrm{x}}{\mathrm{r}}\right\rfloor^2)$

Now $\theta(\mathrm{x})$ can be further written as:

$$\theta(\mathrm{x}) = \sum_{\mathrm{r}\leq\mathrm{x}}(\left(\frac{\mathrm{x}}{\mathrm{r}}\right)^2 + \mathrm{O}(\mathrm{x}/\mathrm{r})) = \frac{\mathrm{x}^2}{2}\left(\zeta(2) + \mathrm{O}\left(\frac{1}{\mathrm{x}}\right)\right) + \mathrm{O}(\mathrm{x}\log\mathrm{x}) \tag{87}$$

as $\mathrm{x}\to\infty$. Hence it can be derived that $R(N)\sim(\pi^2/2).N^2\sim 4\zeta(2)(N^2-\frac{N^2}{4})$ .Hence proved that $\zeta(2) = \pi^2/6$.

**Proof 17:**

Let us now discuss about a proof given by Ivan[11] which is almost similar to that of James D. Harper's[12] simple proof. They used the Fubini theorem for integrals and McLaurin's series expansion for $tanh^{-1}$.

To start with let's first state this basic identity as follows:

$$\frac{1}{2}\log\left(\frac{1+y}{1-y}\right) = \sum_{n=0}^{\infty}\frac{y^{2n+1}}{2n+1}\ for\ |y|<1 \tag{88}$$

Another equality that we shall use is:

$$\int_{-1}^{1}\int_{-1}^{1}\frac{1}{1+2xy+y^2}.dy.dx = \int_{-1}^{1}\int_{-1}^{1}\frac{1}{1+2xy+y^2}.dx.dy \tag{89}$$

Now note that

$$\int_{-1}^{1}\int_{-1}^{1}\frac{1}{1+2xy+y^2}.dy.dx = \int_{-1}^{1}[\frac{\arctan\frac{x+y}{\sqrt{1-x^2}}}{\sqrt{1-x^2}}]_{-1}^{1}.dx = \int_{-1}^{1}\frac{\pi}{2\sqrt{1-x^2}}dx = \frac{\pi^2}{2} \tag{90}$$

The right hand side of the equation can written as follows:

$$\int_{-1}^{1}\int_{-1}^{1}\frac{1}{1+2xy+y^2}.dx.dy = \int_{-1}^{1}[\frac{\log(1+2xy+y^2)}{2y}]_{-1}^{1}.dy = \int_{-1}^{1}\frac{log\frac{1+y}{1-y}}{y}.dy = 2\int_{-1}^{1}\sum_{n=0}^{\infty}\frac{y^{2n+1}}{2n+1} = 4\sum_{n=0}^{\infty}\frac{1}{(2n+1)^2} = 4*\left(\frac{3}{4}\zeta(2)\right) = 3\zeta(2) \tag{91}$$

Therefore it follows that $\zeta(2) = \pi^2/6$.

**Proof 18:**

This is a proof due to Josef Hofbauer[13] .We start off with the simple identity:

$$\frac{1}{sin^2x} = \frac{1}{4\,sin^2(\frac{x}{2})cos^2(\frac{x}{2})} = \frac{1}{4}\left[\frac{1}{sin^2\left(\frac{x}{2}\right)} + \frac{1}{cos^2\left(\frac{x}{2}\right)}\right] = \frac{1}{4}\left[\frac{1}{sin^2\left(\frac{x}{2}\right)} + \frac{1}{sin^2\left(\frac{\pi+x}{2}\right)}\right] \tag{92}$$

Now if we put $x=\pi/2$. The $LHS = \frac{1}{sin^2\pi/2} = 1$ . And the RHS will be

$$\mathrm{RHS} = \left(\frac{1}{4}\right)\left[\frac{1}{sin^2\left(\frac{\pi}{4}\right)} + \frac{1}{sin^2\left(\frac{3\pi}{4}\right)}\right] = \left(\frac{1}{16}\right)\left[\frac{1}{sin^2\left(\frac{\pi}{8}\right)} + \frac{1}{sin^2\left(\frac{3\pi}{8}\right)} + \frac{1}{sin^2\left(\frac{5\pi}{8}\right)} + \frac{1}{sin^2\left(\frac{7\pi}{8}\right)}\right] \tag{93}$$

Now notice that the expression can be expanded into the form mentioned below:

$$\text{RHS} = \frac{1}{4^{\text{n}}}\left[\sum_{\text{k}=0}^{2^{\text{n}}-1}\frac{1}{sin^2\left(\frac{2k+1}{2^{\text{n}+1}}\right)\pi}\right] = \frac{2}{4^{\text{n}}}\left[\sum_{\text{k}=0}^{2^{\text{n}-1}-1}\frac{1}{sin^2\left(\frac{2k+1}{2^{\text{n}+1}}\right)\pi}\right] \quad (94)$$

Taking limit $n \rightarrow \infty$ and using $\lim_{N\rightarrow\infty} N\,sin(x/N) = x$ for $N = 2^n$ and $x = (2k + 1)\pi/2$ yields gives us the following equation:

$$1 = \frac{8}{\pi^2}\sum_{k=0}^{\infty} 1/(2k+1)^2 \quad (95)$$

This as was done earlier many times leads to $\zeta(2) = \pi^2/6$.

## 2. Conclusion

None of these proofs is original; most are well known, though some of them might not be very familiar. Many of the proofs mentioned in the paper were discovered in a survey article by Robin Chapman[14]. Another Proof of the Basel Problem is recently found by author[15-16]. In this paper we discuss some of the notable proofs given by mathematicians to the basal problem. Most of the methods used in this paper are very well known whereas some can be found as proofs of problems present in textbooks.

## 3. References

1. F. Beukers; E. Calabi; J.A.C. Kolk, Sums of generalized harmonic series and volumes, Nieuw Arch. Wisk. No.11(1993) pp. 217-224.
2. Boo Rim Choe; An Elementary Proof of $\sum_{\text{n}=1}^{\infty}\frac{1}{\text{n}^2} = \frac{\pi^2}{6}$ ,The American Mathematical Monthly, Vol. 94, No. 7. (1987), pp. 662-663.
3. T. Apostol, A proof that Euler missed: Evaluating ζ(2) the easy way, Math. Intelligencer 5 (1983) pp.59–60
4. T. Apostol ,Mathematical Analysis (Addison-Wesley, 1974).
5. R. A. Kortram, "Simple proof for $\sum_{\text{k}=1}^{\infty}\frac{1}{\text{k}^2} = \frac{\pi^2}{6}$ and $sin(x) = x\Pi k(1 - x^2/k^2\pi^2)$" Mathematics Magazine, 69, . (1996),pp. 122-125
6. L. Pace, Probabilistically Proving that $\zeta(2) = \frac{\pi^2}{6}$, Amer. Math Monthly, 118 (2011) 641-643
7. Y. Matsuoka, An elementary proof of the formula $\sum\frac{1}{\text{k}^2} = \frac{\pi^2}{6}$, Amer. Math. Monthly 68 (1961) 485– 487
8. E. L. Stark, "Another Proof of the Formula $\sum\frac{1}{\text{k}^2} = \frac{\pi^2}{6}$" The. American Mathematical Monthly, 76, 552-553 (1969).
9. J. M. Borwein and P. B. Borwein, ”Pi and the AGM” , Canadian Mathematical Society Series of Monographs and Advanced Texts, no. 4, John Wiley & Sons, New York, 1998.
10. L.-K. Hua, “Introduction to Number Theory”, Springer,(1982)
11. Mircea Ivan ,“A simple solution to Basel problem”, General Mathematics Vol. 16, No. 4 (2008), pp. 111–113
12. James D. Harper, ”A simple proof of $1 + 1/2^2 + 1/3^2 + \cdots = \frac{\pi^2}{6}$”, The American Mathematical Monthly 109(6) (Jun. - Jul., 2003) 540–541.
13. Josef Hofbauer, “A simple proof of $1 + 1/2^2 + 1/3^2 + \cdots = \frac{\pi^2}{6}$ and related identities”, The American Mathematical Monthly 109(2) (Feb., 2002) 196–200.
14. Robin Chapman, Evaluating ζ(2), Preprint, http: //www.maths.ex.ac.uk/~rjc/etc/zeta2.pdf
15. Ghosh, S. (2021, February 24). Another Proof of Basel Problem. https://doi.org/10.13140/RG.2.2.29833.06249/2
16. Ghosh, Sourangshu. (2020). Solution to the Basel Problem using Calculus of Residues. 10.13140/RG.2.2.29833.06249/2.